\documentclass[12pt]{article}
\usepackage{mathrsfs}

\usepackage{amsmath, epsfig, cite}
\usepackage{amssymb}
\usepackage{amsfonts}
\usepackage{latexsym}
\usepackage{graphicx}

\newtheorem{thm}{Theorem}[section]
\newtheorem{prop}[thm]{Proposition}

\newtheorem{lem}[thm]{Lemma}
\newtheorem{conj}[thm]{Conjecture}

\newtheorem{observation}[thm]{Observation}

\newcommand{\qed}{{\hfill\rule{4pt}{7pt}}}
\def\pf{\noindent {\it Proof.} }

\numberwithin{equation}{section}

\makeatletter \@addtoreset{equation}{section} \makeatother

\begin{document}
\rule{0cm}{1cm}
\begin{center}
{\Large\bf  Solution to a conjecture on the\\[3mm] maximal energy
of bipartite bicyclic graphs\footnote{Supported by NSFC and ``the
Fundamental Research Funds for the Central Universities".}}
\end{center}
 \vskip 2mm \centerline{ Bofeng Huo$^{1,2}$,  Shengjin Ji$^{1}$, \ Xueliang Li$^{1}$, Yongtang Shi$^{1}$ }

\begin{center}
$^{1}$Center for Combinatorics and LPMC\\
Nankai University, Tianjin 300071, China\\

\vskip 2mm Email:huobofeng@mail.nankai.edu.cn;\
jishengjin@mail.nankai.edu.cn;\ lxl@nankai.edu.cn;\
shi@nankai.edu.cn
 \vskip 2mm
 $^{2}$Department of Mathematics and
Information Science \\
Qinghai Normal University, Xining 810008, China
\end{center}

\begin{center}
{\bf Abstract}
\end{center}

{\small The energy of a simple graph $G$, denoted by $E(G)$, is
defined as the sum of the absolute values of all eigenvalues of its
adjacency matrix. Let $C_n$ denote the cycle of order $n$ and
$P^{6,6}_n$ the graph obtained from joining two cycles $C_6$ by a
path $P_{n-12}$ with its two leaves. Let $\mathscr{B}_n$ denote the
class of all bipartite bicyclic graphs but not the graph $R_{a,b}$,
which is obtained from joining two cycles $C_a$ and $C_b$ ($a, b\geq
10$ and $a \equiv b\equiv 2\, (\,\textmd{mod}\, 4)$) by an edge. In
[I. Gutman, D. Vidovi\'{c}, Quest for molecular graphs with maximal
energy: a computer experiment, {\it J. Chem. Inf. Sci.}
{\bf41}(2001), 1002--1005], Gutman and Vidovi\'{c} conjectured that
the bicyclic graph with maximal energy is $P^{6,6}_n$, for $n=14$
and $n\geq 16$. In [X. Li, J. Zhang, On bicyclic graphs with maximal
energy, {\it Linear Algebra Appl.} {\bf427}(2007), 87--98], Li and
Zhang showed that the conjecture is true for graphs in the class
$\mathscr{B}_n$. However, they could not determine which of the two
graphs $R_{a,b}$ and $P^{6,6}_n$ has the maximal value of energy. In
[B. Furtula, S. Radenkovi\'{c}, I. Gutman, Bicyclic molecular graphs
with the greatest energy, {\it J. Serb. Chem. Soc.}
{\bf73(4)}(2008), 431--433], numerical computations up to $a+b=50$
were reported, supporting the conjecture. So, it is still necessary
to have a mathematical proof to this conjecture. This paper is to
show that the energy of $P^{6,6}_n$ is larger than that of
$R_{a,b}$, which proves the conjecture for bipartite bicyclic
graphs. For non-bipartite bicyclic graphs, the conjecture is still
open.}

\vskip 3mm



\section{Introduction}

Let $G$ be a graph of order $n$ and $A(G)$ the adjacency matrix of
$G$. The characteristic polynomial of $G$ is defined as
\begin{equation}\label{i-a}
\phi(G,x)=det(xI-A(G))=\sum_{i=0}^{n}a_{i} x^{n-i}.
\end{equation}
The roots $\lambda_{1}, \lambda_{2}, \ldots, \lambda_{n}$ of
$\phi(G, x)=0$ are called the eigenvalues of $G$.

If $G$ is a bipartite graph, the characteristic polynomial of $G$
has the form
\begin{equation*}
\phi(G,x)=\sum^{\lfloor\frac{n}{2}\rfloor}_{k=0}a_{2k}x^{n-2k}
=\sum^{\lfloor\frac{n}{2}\rfloor}_{k=0}(-1)^{k}b_{2k}x^{n-2k},
\end{equation*}where $b_{2k}=(-1)^{k}a_{2k}$ for all
$k=1,\ldots,\lfloor\frac{n}{2}\rfloor$, especially $b_0=a_0=1$. In
particular, if $G$ is a tree, the characteristic polynomial of $G$
can be expressed as
\begin{equation*}
\phi(G,x)
=\sum^{\lfloor\frac{n}{2}\rfloor}_{k=0}(-1)^{k}m(G,k)x^{n-2k},
\end{equation*}where $m(G,k)$ is the number of $k$-matchings of $G$.

In the following, two basic properties of the characteristic
polynomial $\phi(G)$\cite{C.D1980} will be stated:
\begin{prop}\label{pro1}
If $G_1,G_2,\ldots,G_r$ are the connected components of a graph $G$,
then
\begin{equation*}
\phi(G)=\prod^{r}_{i=1}\phi(G_i).
\end{equation*}
\end{prop}

\begin{prop}\label{pro2}
Let $uv$ be an edge of $G$. Then
\begin{equation*}
\phi(G,x)=\phi(G-uv,x)-\phi(G-u-v,x)-2\sum_{C\in \mathcal
{C}(uv)}\phi(G-C,x),
\end{equation*}where $\mathcal {C}(uv)$ is the set of cycles
containing uv. In particular, if $uv$ is a pendent edge with pendent
vertex $v$, then $\phi(G,x)=x\phi(G-v,x)-\phi(G-u-v,x).$
\end{prop}

The energy of $G$, denoted by $E(G)$, is defined as
$E(G)=\sum\limits^{n}_{i=0}|\lambda_{i}|$. This definition was
proposed by Gutman \cite{gutman1977}. The following formula is also
well-known
\begin{equation*}
E(G)=\frac{1}{\pi}\int^{+\infty}_{-\infty}\frac{1}{x^2}
\log|x^{n}\phi(G,i/x)|\mathrm{d}x,
\end{equation*}where $i^2=-1$. Moreover,
it is known from $[1]$ that the above equality can be expressed as
the following explicit formula:
\begin{equation*}\label{E-dif}
E(G)=\frac{1}{2\pi}\int^{+\infty}_{-\infty}\frac{1}{x^{2}}\log
\left[\left(\sum\limits^{\lfloor n/2\rfloor}_{i=0}(-1)^{i}
a_{2i}x^{2i}\right)^{2}+\left(\sum\limits^{\lfloor n/2\rfloor}_{i=0}(-1)^{i}
a_{2i+1}x^{2i+1}\right)^{2}\right]\mathrm{d}x,
\end{equation*}
where $a_{1}, a_{2}, \ldots, a_{n}$ are the coefficients of the
characteristic polynomial $\phi(G,x)$. For more results about graph
energy, we refer the readers to a survey of Gutman, Li and
Zhang\cite{gutman&li&zhang2009}.

Since 1980s, the extremal energy $E(G)$ of a graph $G$ has been
studied extensively, but the common method makes use of  the
quasi-order. When the graphs are acyclic, bipartite or unicyclic, it
is almost always valid. However, for general graphs, the quasi-order
method is invalid. Recently, for these quasi-order incomparable
problems, we found an efficient way to determine which one attains
the extremal value of the energy, see
\cite{huo&ji&li&shi2011,huo&ji&li2011,huo&ji&li2010,huo&li&shi&wang2011,
huo&li&shi20112u,huo&li&shi2011ubi}, especially, in
\cite{huo&li&shi20112u} we completely solved a conjecture that
$P_n^6$ has the maximal energy among all unicyclic graphs of order
$n\geq 16$.

In this paper, graphs under our consideration are finite, connected
and simple. Let $P_{n}$ and $C_{n}$ denote the path and cycle with
$n$ vertices, respectively. Let $P_{n}^{\ell}$ be the unicyclic
graph obtained by joining a vertex of $C_\ell$ with a leaf of
$P_{n-\ell}$, and $P^{6,6}_n$ the graph obtained from joining two
cycles $C_6$ by a path $P_{n-12}$ with its two leaves. Denote by
$R_{a,b}$ the graph obtained from connecting two cycles $C_a$ and
$C_b$ ($a, b\geq 10$ and $a \equiv b\equiv 2\, (\,\textmd{mod}\,
4)$) by an edge. Let $\mathscr{B}_n$ be the class of all bipartite
bicyclic graphs but not the graph $R_{a,b}$. In
\cite{gutman&vidovic2001}, Gutman and Vidovi\'{c} proposed the
following conjecture on bicyclic graphs with maximal energy:
\begin{conj}\label{conj}
For $n = 14$ and $n\geq 16$, the bicyclic molecular graph of order
$n$ with maximal energy is the molecular graph of the $\alpha,\beta$
diphenyl-polyene $C_6H_5(CH)_{n-12}C_6H_5$, or denoted by
$P^{6,6}_{n}$.
\end{conj}

For bipartite bicyclic graphs, Li and Zhang in \cite{li&zhang2007}
got the following result, giving a partial solution to the above
conjecture.
\begin{thm}\label{mainresultofliz}
If $G\in \mathscr{B}_{n}$, then
$E(G)\leq E(P^{6,6}_{n})$ with equality if and only if $G\cong
P^{6,6}_{n}$.
\end{thm}
However, they could not compare the energies of $P^{6,6}_{n}$ and
$R_{a,b}$. Furtula et al. in \cite{frgutman2008} showed that
$E(P^{6,6}_{n})> E(R_{a,b})$ by numerical computations up to
$a+b=50$, supporting that the conjecture is true for bipartite
bicyclic graphs. It is evident that a mathematical proof is still
needed. This paper is to give such a proof. We will use Coulson
integral formula and some knowledge of real analysis as well as
combinatorial method to show the following result:
\begin{thm}\label{mainresultmy}
For $n-t,t\geq 10$ and $n-t\equiv t\equiv 2\, \ (mod\, 4)$,
$E(R_{n-t,t})<E(P^{6,6}_{n})$.
\end{thm}
As Furtula et al. noticed in \cite{frgutman2008}, since for odd $n$
the graph $R_{a,b} \ (a+b=n)$ is not bipartite, therefore, for odd
$n$, it is known that $P^{6,6}_{n}$ is the maximal energy bipartite
bicyclic graph from \cite{li&zhang2007}. Therefore, combining
Theorems \ref{mainresultofliz} and \ref{mainresultmy}, we get:
\begin{thm}
Let $G$ be any connected, bipartite bicyclic graph with $n\,(\,n\geq
12)$ vertices. Then $E(G)\leq E(P^{6,6}_{n})$ with equality if and
only if $G\cong P^{6,6}_{n}$.
\end{thm}

So, Conjecture \ref{conj} is true for all connected bipartite
bicyclic graphs of order $n$ with $n = 14$ and $n\geq 16$. However,
it is still open for non-bipartite bicyclic graphs.

\section{Proof of Theorem \ref{mainresultmy}}

Before giving the proof of Theorem \ref{mainresultmy}, we shall
state some knowledge on real analysis \cite{Zorich}.
\begin{lem}\label{X1}
For any real number $X>-1$, we have
\begin{equation*}\label{X}
\frac{X}{1+X}\leq \log(1+X)\leq X.
\end{equation*}In particular, $\log(1+X)<0$ if and only if $X<0$.
\end{lem}
The following lemma is a well-known conclusion due to Gutman
\cite{gutman2001} which will be used later.
\begin{lem}\label{E-E}
If $G_1$ and $G_2$ are two graphs with the same number of vertices,
then
\begin{equation*}
E(G_1)-E(G_2)=\frac{1}{\pi}\int^{+\infty}_{-\infty}
\log\frac{\phi(G_1;ix)}{\phi(G_2;ix)}\mathrm{d}x.
\end{equation*}
\end{lem}

One can easily obtain the following recursive equations from
Propositions \ref{pro1} and \ref{pro2}.
\begin{lem}\label{polynsome}
For any positive number $n\geq 8$,
\begin{alignat*}{1}
\phi(P_n,x)&=x\phi(P_{n-1},x)-\phi(P_{n-2},x),\\
\phi(C_n,x)&=\phi(P_{n},x)-\phi(P_{n-2},x)-2,\\
\phi(P^6_n,x)&=x\phi(P^6_{n-1},x)-\phi(P^6_{n-2},x);
\end{alignat*}
for any positive number $n\geq 6$ and $t\geq 3$,
\begin{equation*}
  \phi(R_{n-t,t},x)=\phi(C_{n-t},x)\phi(C_t,x)-\phi(P_{n-t-1},x)\phi(P_{t-1},x).
\end{equation*}

\end{lem}

Next, we introduce some convenient notations as follows, which will
be used in the sequel.
\begin{equation*}
Y_1(x)=\frac{x+\sqrt{x^2-4}}{2},\qquad\qquad
Y_2(x)=\frac{x-\sqrt{x^2-4}}{2}.
\end{equation*}It is easy to verify that $Y_1(x)+Y_2(x)=x$,
$Y_1(x)Y_2(x)=1$, $Y_1(ix)=\frac{x+\sqrt{x^2+4}}{2}i$ and
$Y_2(ix)=\frac{x-\sqrt{x^2+4}}{2}i$. Furthermore, we define
\begin{equation*}
Z_1(x)=-iY_1(ix)=\frac{x+\sqrt{x^2+4}}{2},\quad
Z_2(x)=-iY_2(ix)=\frac{x-\sqrt{x^2+4}}{2}.
\end{equation*}Note that $Z_1(x)+Z_2(x)=x$,
$Z_1(x)Z_2(x)=-1$. Moreover, $Z_1(x)>1$ and $-1<Z_2(x)<0$, if $x>0$;
$0<Z_1(x)<1$ and $Z_2(x)<-1$, otherwise. In the rest of this paper,
we abbreviate $Z_{j}(x)$ to $Z_{j}$ for $j=1,2.$ Some more notations
will be used frequently in the sequel.
\begin{alignat*}{2}
A_1(x)&=\frac{Y_1(x)\phi(P^{6,6}_{13},x)-\phi(P^{6,6}_{12},x)}{(Y_1(x))^{14}-(Y_1(x))^{12}},&
\quad
A_2(x)&=\frac{Y_2(x)\phi(P^{6,6}_{13},x)-\phi(P^{6,6}_{12},x)}{(Y_2(x))^{14}-(Y_2(x))^{12}},\\
B_1(x)&=\frac{Y_1(x)(x^2-1)-x}{(Y_1(x))^3-Y_1(x)},&
B_2(x)&=\frac{Y_2(x)(x^2-1)-x}{(Y_2(x))^3-Y_2(x)}.
\end{alignat*}By some simple calculations, we have that
$\phi(P^{6,6}_{13},x)=x^{13}-14x^{11}+74x^9-188x^7+245x^5-158x^3+40x$
and $\phi(P^{6,6}_{12},x)=x^{12}-13x^{10}+62x^8-138x^6+153x^4-81x^2+16$,
and then
\begin{equation*}
A_1(ix)=\frac{Z_{1}g_{13}+g_{12}}{Z_1^2+1}Z_2^{12},\qquad
A_2(ix)=\frac{Z_{2}g_{13}+g_{12}}{Z_1^2+1}Z_1^{12},
\end{equation*}where
$g_{13}=x^{13}+14x^{11}+74x^9+188x^7+245x^5+158x^3+40x$ and
$g_{12}=x^{12}+13x^{10}+62x^8+138x^6+153x^4+81x^2+16$. Notice that
$A_j(ix)$ has a good property, i.e., its sign is always positive for
all real number $x$, for $j=1,2$.
\begin{observation}
For all real number $x$, $A_j(ix)>0$, $j=1,2.$
\end{observation}
\pf Since, by some directed calculations, we have
\begin{equation*}
A_1(ix)A_2(ix)=\frac{(x^6+8x^4+19x^2+16)^2(x^2+1)^4}{x^2+4}>0 \
\text{for\, all}\  x.
\end{equation*}
Besides, from the expression of $A_1(ix)$, we obviously obtain that
$A_1(ix)>0$ for all real $x$. Thus, we conclude that $A_2(ix)>0$.
For convenience, we abbreviate $A_{j}(ix)$ and $C_{j}(ix)$ to
$A_{j}$ and $C_{j}$ for $j=1,2$, respectively.\qed

 The following lemma will be used in the showing of the later results, due to Huo et al.
\cite{huo&ji&li&shi2011,huo&li&shi2011ubi,huo&li&shi20112u}.
\begin{lem}\label{polypc}
For $n\geq 4$ and $x\neq \pm 2$, the characteristic polynomials of
$P_n$ and $C_n$ possess  the following forms:
\begin{equation*}
\phi(P_n,x)=B_1(x)(Y_1(x))^n+B_2(x)(Y_2(x))^n
\end{equation*}
and
\begin{equation*} \phi(C_n,x)=(Y_1(x))^n+(Y_2(x))^n-2.
\end{equation*}
\end{lem}
\begin{lem}
For $n\geq 12$, the characteristic polynomial of $P^{6,6}_n$ has the
following form:
\begin{equation*}
\phi(P^{6,6}_n,x)=A_1(x)(Y_1(x))^n+A_2(x)(Y_2(x))^n
\end{equation*}where $x\neq \pm 2$.
\end{lem}
\pf Note that, $\phi(P^{6,6}_n)$ satisfies the recursive formula
$f(n,x)=xf(n-1,x)-f(n-2,x)$ in terms of Lemma \ref{polynsome}.
Therefore, the form of the general solution of the linear
homogeneous recursive relation is
$f(n,x)=D_1(x)(Y_1(x))^n+D_2(x)(Y_2(x))^n$. By some simple
calculations, together with the initial values $\phi(P^{6,6}_{12})$
and $\phi(P^{6,6}_{13})$, we can get that $D_i(x)=A_i(x)$, $i=1,2.$
\qed

From Lemmas \ref{polynsome} and \ref{polypc} and Proposition
\ref{pro1}, by means of elementary calculations it is easy to deduce
the following result. The details of its proof is omitted.
\begin{lem}\label{polynrnt}
For $n\geq 6$ and $t\geq 3$, the characteristic polynomial of
$R_{n-t,t}$ has the following form:
\begin{equation*}
\phi(R_{n-t,t},x)=C_1(x)(Y_1(x))^n+C_2(x)(Y_2(x))^n-2((Y_1(x))^t+(Y_2(x))^t)+4
\end{equation*}where $x\neq \pm 2$,
$C_1(x)=1+(Y_2(x))^{2t}-2(Y_2(x))^t-(B_1(x))^2(Y_2(x))^2-B_1(x)B_2(x)(Y_2(x))^{2t}$
and
$C_2(x)=1+(Y_1(x))^{2t}-2(Y_1(x))^t-(B_2(x))^2(Y_1(x))^2-B_1(x)B_2(x)(Y_1(x))^{2t}$.
\end{lem}
In terms of the above lemma, we can get the following forms for
$C_j(ix)$\,($j=1,2$) by some simplifications,
\begin{alignat*}{1}
C_1(ix)&=1+\frac{x^2+3}{x^2+4}Z_2^{2t}+2Z_2^t+\frac{Z_1^2}{(Z_1^2+1)^2}\\
C_2(ix)&=1+\frac{x^2+3}{x^2+4}Z_1^{2t}+2Z_1^t+\frac{Z_2^2}{(Z_2^2+1)^2}.
\end{alignat*}
\\

{\bf Proof of Theorem \ref{mainresultmy}}\\

From the above analysis, we only need to show that
$E(R_{n-t,t})<E(P^{6,6}_n)$, for every positive number $t=4k_1+2 \
(t\geq 10)$, $n-t\geq 10$ and $n=4k_2 \ (n\geq 2t)$. Without loss of
generality, we assume $n-t\geq t$, that is, $n\geq 2t$. From Lemma
\ref{E-E}, we have that
\begin{equation*}
E(R_{n-t,t})-E(P^{6,6}_n)=\frac{1}{\pi}\int^{+\infty}_{-\infty}
\log\frac{\phi(R_{n-t,t};ix)}{\phi(P^{6,6}_n;ix)}\mathrm{d}x.
\end{equation*}
First of all, we shall will that the integrand
$\log\frac{\phi(R_{n-t,t};ix)}{\phi(P^{6,6}_n;ix)}$ is monotonically
decreasing in $n$ for $n=4k$, that is,
\begin{alignat*}{1}
 &\log\frac{\phi(R_{n+4-t,t};ix)}{\phi(P^{6,6}_{n+4};ix)}-
  \log\frac{\phi(R_{n-t,t};ix)}{\phi(P^{6,6}_n;ix)} \\
  &=\log\frac{\phi(R_{n+4-t,t};ix)\phi(P^{6,6}_n;ix)}
  {\phi(P^{6,6}_{n+4};ix)\phi(R_{n-t,t};ix)}
  =\log\left(1+\frac{K(n,t,x)}{H(n,t,x)}\right),
\end{alignat*}where
$K(n,t,x)=\phi(R_{n+4-t,t};ix)\phi(P^{6,6}_n;ix)-\phi(P^{6,6}_{n+4};ix)\phi(R_{n-t,t};ix)$
and $H(n,t,x)=\phi(P^{6,6}_{n+4};ix)\phi(R_{n-t,t};ix)>0.$ From
Lemma \ref{X}, we only need to verify that $K(n,t,x)<0$. By means of
some directed calculations, we arrive at
\begin{equation*}
K(n,t,x)=(Z_1^4-Z_2^4)(A_2C_1-A_1C_2)+(2Z_1^t+2Z_2^t+4)(A_1Z_1^n(1-Z_1^4)+A_2Z_2^n(1-Z_2^4)).
\end{equation*}
Noticing that $Z_1>1$ and $0>Z_2>-1$ for $x>0$, we have $Z_1^n\geq
Z_1^{2t}>0$ and $0<Z_2^n\leq \linebreak Z_2^{2t}$. Meanwhile, from
$0<Z_1<1$ and $Z_2<-1$ for $x<0$, we have $0<Z_1^n\leq Z_1^{2t}$ and
$Z_2^n\geq Z_2^{2t}>0$. Therefore,
\begin{equation*}
A_1Z_1^n(1-Z_1^4)+A_2Z_2^n(1-Z_2^4)\leq
A_1Z_1^{2t}(1-Z_1^4)+A_2Z_2^{2t}(1-Z_2^4).
\end{equation*}Namely, $K(n,t,x)\leq K(2t,t,x)=(Z_1^4-Z_2^4)(A_2C_1-A_1C_2)+(2Z_1^t+2Z_2^t+4)
(A_1Z_1^{2t}(1-Z_1^4)+A_2Z_2^{2t}(1-Z_2^4)).$ Now let
$f(t,x)=K(2t,t,x)$. By some simplifications, it is easy to get
\begin{equation*}
  f(t,x)=\alpha_0Z_1^{3t}+\alpha_1Z_1^{-3t}+\beta_0Z_1^{2t}+\beta_1Z_1^{-2t}
  +\gamma_0Z_1^{t}+\gamma_1Z_1^{-t}+a_0,
\end{equation*}where
\begin{alignat*}{2}
\alpha_0&=2A_1(1-Z_1^{4}),&\quad
\alpha_1&=2A_2(1-Z_2^{4}), \\
\beta_0&=A_1\left(4(1-Z_1^{4})-(Z_1^{4}-Z_2^{4})\frac{x^2+3}{x^2+4}\right),&
\beta_1&=A_2\left(4(1-Z_2^{4})+(Z_1^{4}-Z_2^{4})\frac{x^2+3}{x^2+4}\right), \\
\gamma_0&=2A_1((1-Z_1^{4})-(Z_1^{4}-Z_2^{4})),&
\gamma_1&=2A_2((1-Z_2^{4})+(Z_1^{4}-Z_2^{4})),
\end{alignat*}and
\begin{alignat*}{2}
  a_0&=(Z_1^{4}-Z_2^{4})\left(A_2\left(1+\frac{Z_1^{2}}{(Z_1^{2}+1)^2}\right)
-A_1\left(1+\frac{Z_2^{2}}{(Z_2^{2}+1)^2}\right)\right).&\qquad\qquad\qquad\qquad\quad\
&
\end{alignat*} {\bf Claim 1.} $f(t,x)$ is monotonically decreasing in $t$.

Note the facts that $(1-Z_1^{4})<0$ for $x>0$, $(1-Z_1^{4})>0$ for
$x<0$; $(1-Z_2^{4})>0$ for $x>0$, $(1-Z_2^{4})<0$ for $x<0$;
$(Z_1^{4}-Z_2^{4})>0$ for $x>0$, $(Z_1^{4}-Z_2^{4})<0$ for $x<0$. It
is not difficult to check that $\alpha_0<0, \beta_0<0 \
\textmd{and}\ \gamma_0<0$ for $x>0$, $\alpha_0>0, \beta_0>0 \
\textmd{and}\ \gamma_0>0$, otherwise; thus $\alpha_1>0, \beta_1>0 \
\textmd{and}\ \gamma_1>0$ for $x>0$, $\alpha_1<0, \beta_1<0 \
\textmd{and}\ \gamma_1<0$, otherwise. Therefore, no matter which of
$x>0$ or $x<0$ happens, we can always deduce that
 \begin{equation*}
 \frac{\partial f(t,x)}{\partial t}=(3\alpha_0Z_1^{3t}-3\alpha_1Z_1^{-3t}+2\beta_0Z_1^{2t}-2\beta_0Z_1^{-2t}
  +\gamma_0Z_1^{t}-\gamma_1Z_1^{-t})\log Z_1<0.
\end{equation*}
The proof of Claim\,1 is thus complete.

From Claim 1, it follows that for $t\geq 10$,
\begin{alignat*}{1}
K(n,t,x)\leq f(10,x)=&-4x^2(x^2+1)^2(x^{18}+23x^{16}+224x^{14}+1203x^{12}+3887x^{10}\\
                      & +7731x^8+9285x^6+6301x^4+2077x^2+224)\\
                      &-(x^{10}+13x^8+62x^6+131x^4+109x^2+16)\\
                      & \times
                      x^2(x^4+5x^2+6)(x^4+3x^2+1)(x^2+1)^2<0.
\end{alignat*} Therefore, we have verified that the integrand
$\log\frac{\phi(R_{n-t,t};ix)}{\phi(P^{6,6}_n;ix)}$ is monotonically
decreasing in $n$ for $n=4k$. That is,
$E(R_{n-t,t})-E(P^{6,6}_n)\leq E(R_{10,10})-E(P^{6,6}_{12})<0$ for
every positive number $t=4k_1+2 \ (n\geq 10)$, $n-t\geq 10$ and
$n=4k_2 \ (n\geq 2t)$. Therefore, the entire proof of
Theorem\ref{mainresultmy} is now complete.\qed

\noindent{\bf Acknowledgement.} The authors are very grateful to the
referees for helpful comments and suggestions, which improved the
presentation of the original manuscript.

\end{document}